\documentclass[a4paper, 11pt]{article}

\usepackage{amsmath}
\usepackage{amsfonts}
\usepackage{amssymb}
\usepackage[english]{babel}
\usepackage{graphicx}
\usepackage{amsthm}

\theoremstyle{definition}

\begin{document}

\title{Octahedral, dicyclic and special linear solutions \break of some unsolved Hamilton-Waterloo problems}

\author{S.Bonvicini, \thanks{Dipartimento di Scienze e Metodi dell'Ingegneria,
Universit\`a di Modena e Reggio Emilia, via Amendola 2, 42100 Reggio Emilia, Italy, email: simona.bonvicini@unimore.it} \quad
M. Buratti \thanks{Dipartimento di Matematica e Informatica, Universit\`{a} degli Studi di Perugia, Via Vanvitelli 1, 06123 Perugia, Italy. email: buratti@dmi.unipg.it}}

\maketitle

\begin{abstract}
\noindent We give a sharply-vertex-transitive solution of each of the nine Hamilton-Waterloo problems
left open by Danziger, Quattrocchi and Stevens.
\end{abstract}

\section{Introduction}

In a long paper still in preparation \cite{BurBon}, we will give a general method to obtain decompositions and/or factorizations of a graph with a nice automorphism group, in
particular with an automorphism group $G$ acting sharply transitively (i.e., regularly) on the vertices. In this case one also say that the decomposition or factorization is
$G$-regular.

Sharply-vertex-transitive 2-factorizations of a complete graph of odd order have been already partly investigated in \cite{BR}.
Sharply-vertex-transitive 2-factorizations of a cocktail party graph (that is a complete graph of even order minus
a 1-factor) can be treated similarly. Here we show how the method of {\it partial differences} explained in \cite{B2004} and
successfully applied to solve  several cycle decompositions problems (see, e.g., \cite{B}) allows to obtain a
sharply-vertex-transitive solution of each of the nine Hamilton-Waterloo problems left open by Danziger, Quattrocchi and Stevens
\cite{DQS}. We namely prove that: 1) there exists a regular 2-factorization of $K_{48}-I$ having exactly $r$ triangle-factors
with $r\in\{5,7,9, 13,15,17\}$ and each of the remaining factors consisting of all quadrangles; 2) there exists a regular
2-factorization of $K_{24}-I$ having exactly $r$ triangle-factors  with $r\in\{5,7,9\}$ and each of the remaining
factors consisting of all quadrangles. The choice of the group acting regularly on our factorizations cannot be random;
in each case we have been forced to take an ad hoc non-abelian group with exactly one involution.

\section {Octahedral solutions of six Hamilton-Waterloo problems}

Throughout this section $G$ will denote the so-called {\it binary octahedral group} which is usually denoted by $\overline O$.
This group, up to isomorphism, can be viewed as a group of units
of the skew-field $\mathbb{H}$ of \emph{quaternions} introduced by Hamilton, that is an extension of
the complex field $\mathbb{C}$. We recall the basic facts regarding
$\mathbb{H}$. Its elements are all real linear combinations of $1$,
$i$, $j$ and $k$. The sum and the product of two quaternions are
defined in the natural way under the rules that
$$i^2=j^2=k^2=ijk=-1.$$
If $q=a+bi+cj+dk \neq 0$, its inverse is given by $$q^{-1}={a-bi-cj-dk\over a^2+b^2+c^2+d^2}.$$
The elements of the
multiplicative group $G$ can be described as follows:

$$
\begin{array}{l}
\pm 1, \pm i, \pm j, \pm k;\smallskip\\
\frac 12(\pm 1\pm i\pm j\pm k);\smallskip\\
\frac 1{\sqrt 2}(\pm x\pm y), \quad \{x, y\}\in\binom{\{1, i, j, k\}}2.\\
\end{array}
$$

The use of octahedral group $G$ was crucial in \cite{BBRT} to get a Steiner triple system of any order $v=96n+49$
with an automorphism group acting sharply transitively an all but one point.
Here $G$ will be used to get a $G$-regular solution of each of the six Hamilton-Waterloo problems of order 48 left open in \cite{DQS}.
We will need to consider the following subgroups of $G$ order $16$ and $12$,
respectively:
\begin{itemize}
\item $K=\langle k, \  {1\over \sqrt{2}}(j-k)\rangle$;
\item $L=\langle {1\over \sqrt{2}}(j-k), \ \frac 12(-1-i+j+k) \rangle$.
\end{itemize}

In the following the complete graph $K_{48}$ and the cocktail party graph $K_{48}-I$
will be seen as the Cayley graphs Cay$[G:G\setminus\{1\}]$ and\break Cay$[G:G\setminus\{1,-1\}]$, respectively.

\subsection{An octahedral solution of HWP$(48; 3, 4; 5, 18)$}

\noindent
Consider the nine cycles of $K_{48}$ defined as follows.

$$
\begin{array}{ll}
C_1=&\bigl{(}1, \ -\frac 1{\sqrt 2}(1-k), \ \frac 1{2}(1-i-j-k)\bigl{)}\smallskip\\
C_2=&\bigl{(}1, \ \frac 12 (-1-i+j+k), \ \frac 12 (-1+i-j-k)\bigl{)}\smallskip\\
C_3=&\bigl{(}1, \ \frac 12 (-1+i+j-k), \ \frac 12 (-1-i-j+k)\bigl{)}\smallskip\\
C_4=&\bigl{(}1, \ k, \ -1, \ -k\bigl{)}\smallskip\\
C_5=&\bigl{(}1, \ j, \ -1, \ -j\bigl{)}\smallskip\\
C_6=&\bigl{(}1, \ \frac 1{\sqrt 2}(-i+k), \ -\frac 1{2}(1+i+j+k), \ -\frac 1{\sqrt 2}(j+k)\bigl{)}\smallskip\\
C_7=&\bigl{(}1, \ \frac 1{\sqrt 2}(i-j), \ \frac 1{\sqrt 2}(1+i), \ \frac 1{2}(1-i-j+k)\bigl{)}\smallskip\\
C_8=&\bigl{(}1, \ \frac 1{2}(1-i+j-k), \ k, \ -\frac 1{\sqrt 2}(1+j)\bigl{)}\smallskip\\
C_9=&\bigl{(}1, \ \frac 1{\sqrt 2}(1-i), \ -\frac 1{\sqrt 2}(1+i), \ \frac 1{2}(-1-i+j-k)\bigl{)}\smallskip\\
\end{array}
$$

We note that for $2\leq i\leq 5$ the $G$-stabilizer of $C_i$ is $V(C_i)$ while all other $C_i$'s have trivial $G$-stabilizer.
Thus, using partial differences, one can check that $Orb_G(C_i)$ is a $\ell_i$-cycle decomposition of Cay$[G:\Omega_i]$ where $\ell_i$
is the length of $C_i$ and where the $\Omega_i$'s are the symmetric subsets of $G$ listed below.

$$
\begin{array}{ll}
\Omega_1=&\{-\frac 1{\sqrt 2}(1-k), \frac 1{2}(1-i-j-k), -\frac 1{\sqrt 2}(1+i)\}^{\pm 1}\smallskip\\
\Omega_2=&\{\frac 12 (-1-i+j+k)\}^{\pm 1}\smallskip\\
\Omega_3=&\{\frac 12 (-1+i+j-k)\}^{\pm 1}\smallskip\\
\Omega_4=&\{k\}^{\pm 1}\smallskip\\
\Omega_5=&\{j\}^{\pm 1}\smallskip\\
\Omega_6=&\{\frac 1{\sqrt 2}(-i+k), \frac 1{\sqrt 2}(j-k), \frac 1{\sqrt 2}(1-k), -\frac 1{\sqrt 2}(j+k)\}^{\pm 1}\smallskip\\
\Omega_7=&\{\frac 1{\sqrt 2}(i-j), \frac 12(1+i-j-k), \frac 1{\sqrt 2}(i+j), \frac 1{2}(1-i-j+k)\}^{\pm 1}\smallskip\\
\Omega_8=&\{\frac 1{2}(1-i+j-k), -\frac 1{2}(1+i+j+k), -\frac 1{\sqrt 2}(i+k), -\frac 1{\sqrt 2}(1+j)\}^{\pm 1}\smallskip\\
\Omega_9=&\{\frac 1{\sqrt 2}(1-i), i, \frac 1{\sqrt 2}(1+j), \frac 1{2}(-1-i+j-k)\}^{\pm 1}\smallskip\\
\end{array}
$$

Now note that the $\Omega_i$'s partition $G\smallsetminus\{1, -1\}$.
This assures that \break $\displaystyle {\cal C}:= \bigcup_{i=1}^{7}Orb_{G}(C_i)$ is a $G$-regular
cycle-decomposition of $K_{48}-I$.

Now set $F_i=Orb_{S_i}(C_i)$ where
$S_i=\begin{cases}K & \mbox{for $i=1$;}\cr G & \mbox{for $2\leq i \leq 5$;}\cr L & \mbox{for $6 \leq i\leq 9$}.\end{cases}$

Each $F_i$ is a 2-factor of $K_{48}$ with $Stab_{G}(F_i)=S_i$, hence $Orb_{G}(F_i)$ has length
3 or 1 or 4 according to whether $i=1$, or $2\leq i\leq 5$, or $6\leq i\leq 9$, respectively.
The cycles of $F_i$ are triangles or quadrangles according to whether or not $i\leq3$.
Thus, recalling that $\cal C$ is a cycle-decomposition of $K_{48}-I$, we conclude that
$\displaystyle {\cal F}:= \bigcup_{i=1}^{9}Orb_{G}(F_i)$ is a $G$-regular
2-factorization of $K_{48}-I$ with 5 triangle-factors and 18 quadrangle-factors, namely a
$G$-regular solution of HWP$(48; 3, 4; 5, 18)$.

\subsection{An octahedral solution of HWP$(48; 3, 4; 7, 16)$}

\noindent
Consider the seven cycles of $K_{48}$ defined as follows.

$$
\begin{array}{ll}
C_1=&\bigl{(}1, \ -\frac 1{\sqrt 2}(i+j), \ \frac 1{2}(1-i+j+k)\bigl{)}\smallskip\\
C_2=&\bigl{(}1, \ \frac 1{2}(-1-i+j+k), \ \frac 1{2}(1-i-j-k)\bigl{)}\smallskip\\
C_3=&\bigl{(}1, \ \frac 12 (-1+i+j-k), \ \frac 12 (-1-i-j+k)\bigl{)}\smallskip\\
C_4=&\bigl{(}1, \ \frac 1{\sqrt 2}(-i+k), \ \frac 1{2}(1+i+j-k), \ -\frac 1{\sqrt 2}(j+k)\bigl{)}\smallskip\\
C_5=&\bigl{(}1, \ \frac 1{\sqrt 2}(i-j), \ \frac 1{\sqrt 2}(1-k), \ \frac 1{\sqrt 2}(1+i)\bigl{)}\smallskip\\
C_6=&\bigl{(}1, \ \frac 1{\sqrt 2}(1+k), \ -\frac 1{2}(1+i+j+k), \ \frac 1{\sqrt 2}(1+j)\bigl{)}\smallskip\\
C_7=&\bigl{(}1, \ -\frac 1{2}(1+i+j+k), \ \frac 1{2}(1-i+j-k), \ \frac 1{2}(1-i-j+k)\bigl{)}\smallskip\\
\end{array}
$$

We note that the $G$-stabilizer of $C_3$ is $V(C_3)$ while all other $C_i$'s have trivial $G$-stabilizer.
Thus, using partial differences, one can check that $Orb_G(C_i)$ is a $\ell_i$-cycle decomposition of Cay$[G:\Omega_i]$ where $\ell_i$
is the length of $C_i$ and where the $\Omega_i$'s
are the symmetric subsets of $G$ listed below.

$$
\begin{array}{ll}
\Omega_1=&\{-\frac 1{\sqrt 2}(i+j), \frac 1{2}(1-i+j+k), \frac 1{\sqrt 2}(-j+k)\}^{\pm 1}\smallskip\\
\Omega_2=&\{\frac 1{2}(-1-i+j+k), \frac 1{2}(1-i-j-k), \frac 1{2}(-1-i+j-k)\}^{\pm 1}\smallskip\\
\Omega_3=&\{\frac 12 (-1+i+j-k)\}^{\pm 1}\smallskip\\
\Omega_4=&\{\frac 1{\sqrt 2}(-i+k), -\frac 1{\sqrt 2}(1-k), \frac 1{\sqrt 2}(i+k), -\frac 1{\sqrt 2}(j+k)\}^{\pm 1}\smallskip\\
\Omega_5=&\{\frac 1{\sqrt 2}(i-j), -j, \frac 1{2}(1-i+j-k), \frac 1{\sqrt 2}(1+i)\}^{\pm 1}\smallskip\\
\Omega_6=&\{\frac 1{\sqrt 2}(1+k), \frac 1{\sqrt 2}(-1+j), -\frac 1{\sqrt 2}(1+i), \frac 1{\sqrt 2}(1+j)\}^{\pm 1}\smallskip\\
\Omega_7=&\{-\frac 1{2}(1+i+j+k), -i, -k, \frac 1{2}(1-i-j+k)\}^{\pm 1}\smallskip\\
\end{array}
$$

One can see that the $\Omega_i$'s partition $G\smallsetminus\{1, -1\}$.
This assures that $\displaystyle {\cal C}:= \bigcup_{i=1}^{7}Orb_{G}(C_i)$ is an $G$-regular
cycle-decomposition of $K_{48}-I$.

Now set $F_i=Orb_{S_i}(C_i)$ where
$S_i=\begin{cases}K & \mbox{for $i=1, 2$;}\cr G & \mbox{for $i=3$;}\cr L & \mbox{for $4 \leq i\leq 7$}.\end{cases}$

Each $F_i$ is a 2-factor of $K_{48}$ with $Stab_{G}(F_i)=S_i$, hence $Orb_{G}(F_i)$ has length
3 or 1 or 4 according to whether $i=1, 2$ or $i=3$ or $4\leq i\leq 7$, respectively.

The cycles of $F_i$ are triangles or quadrangles according to whether or not $i\leq3$.
Thus, recalling that $\cal C$ is a cycle-decomposition of $K_{48}-I$, we conclude that
$\displaystyle {\cal F}:= \bigcup_{i=1}^{7}Orb_{G}(F_i)$ is a $G$-regular
2-factorization of $K_{48}-I$ with 7 triangle-factors and 16 quadrangle-factors, namely a
$G$-regular solution of HWP$(48; 3, 4; 7, 16)$.

\subsection{An octahedral solution of HWP$(48; 3, 4; 9, 14)$}

\noindent
Consider the eight cycles of $K_{48}$ defined as follows.
$$
\begin{array}{ll}
C_1=&\bigl{(}1, \ \frac 1{\sqrt 2}(i+j), \ \frac 1{2}(1-i-j-k)\bigl{)}\smallskip\\
C_2=&\bigl{(}1, \ -\frac 1{\sqrt 2}(1-k), \ \frac 1{\sqrt 2}(1+j)\bigl{)}\smallskip\\
C_3=&\bigl{(}1, \ \frac 1{2}(-1-i+j+k), \ \frac 1{2}(1+i-j+k)\bigl{)}\smallskip\\
C_4=&\bigl{(}1, \ \frac 1{\sqrt 2}(-i+k), \ \frac 1{\sqrt 2}(1-i), \ \frac 1{2}(-1-i+j-k)\bigl{)}\smallskip\\
C_5=&\bigl{(}1, \ \frac 1{\sqrt 2}(i-j), \ \frac 1{2}(-1+i+j+k), \ -\frac 1{\sqrt 2}(j+k)\bigl{)}\smallskip\\
C_6=&\bigl{(}1, \ \frac 1{\sqrt 2}(1+i), \ \frac 1{\sqrt 2}(1-i), \ \frac 1{2}(1-i-j+k)\bigl{)}\smallskip\\
C_7=&\bigl{(}1, \ k, \ -1, \ -k\bigl{)}\smallskip\\
C_8=&\bigl{(}1, \ j, \ -1, \ -j\bigl{)}\smallskip\\
\end{array}
$$

We note that for $i=7, 8$ the $G$-stabilizer of $C_i$ is $V(C_i)$ while all other $C_i$'s have trivial $G$-stabilizer.
Thus, using partial differences, one can check that $Orb_G(C_i)$ is a $\ell_i$-cycle
decomposition of Cay$[G:\Omega_i]$ where $\ell_i$
is the length of $C_i$ and where the $\Omega_i$'s
are the symmetric subsets of $G$ listed below.

$$
\begin{array}{ll}
\Omega_1=&\{\frac 1{\sqrt 2}(i+j), \frac 1{2}(1-i-j-k), \frac 1{\sqrt 2}(-1+i)\}^{\pm 1}\smallskip\\
\Omega_2=&\{-\frac 1{\sqrt 2}(1-k), \frac 1{\sqrt 2}(1+j), \frac 1{2}(-1+i+j+k)\}^{\pm 1}\smallskip\\
\Omega_3=&\{\frac 1{2}(-1-i+j+k), \frac 1{2}(1+i-j+k), \frac 1{2}(-1-i-j+k)\}^{\pm 1}\smallskip\\
\Omega_4=&\{\frac 1{\sqrt 2}(-i+k), \frac 1{2}(1-i+j+k), \frac 1{\sqrt 2}(i+k), \frac 12(-1-i+j-k)\}^{\pm 1}\smallskip\\
\Omega_5=&\{\frac 1{\sqrt 2}(i-j), \frac 1{\sqrt 2}(j-k), -\frac 1{\sqrt 2}(1+j), -\frac 1{\sqrt 2}(j+k)\}^{\pm 1}\smallskip\\
\Omega_6=&\{\frac 1{\sqrt 2}(1+i), i, \frac 1{\sqrt 2}(1-k), \frac 1{2}(1-i-j+k)\}^{\pm 1}\smallskip\\
\Omega_7=&\{k\}^{\pm 1}\smallskip\\
\Omega_8=&\{j\}^{\pm 1}\smallskip\\
\end{array}
$$

Now note that the $\Omega_i$'s partition $G\smallsetminus\{1, -1\}$.
This assures that \break $\displaystyle {\cal C}:= \bigcup_{i=1}^{8}Orb_{G}(C_i)$ is a $G$-regular
cycle-decomposition of $K_{48}-I$.

Now set $F_i=Orb_{S_i}(C_i)$ where
$S_i=\begin{cases}K & \mbox{for $1 \leq i \leq 3$;}\cr L & \mbox{for $4\leq i\leq 6$;}\cr G & \mbox{for $i=7, 8$}.\end{cases}$

Each $F_i$ is a 2-factor of $K_{48}$ with $Stab_{G}(F_i)=S_i$, hence $Orb_{G}(F_i)$ has length
3 or 4 or 1 according to whether $1\leq i\leq 3$ or $4 \leq i \leq 6$ or $i=7, 8$, respectively.
The cycles of $F_i$ are triangles or quadrangles according to whether or not $i\leq3$.
Thus, recalling that $\cal C$ is a cycle-decomposition of $K_{48}-I$, we conclude that
$\displaystyle {\cal F}:= \bigcup_{i=1}^{8}Orb_{G}(F_i)$ is a $G$-regular
2-factorization of $K_{48}-I$ with 9 triangle-factors and 14 quadrangle-factors, namely a
$G$-regular solution of HWP$(48; 3, 4; 9, 14)$.

\subsection{An octahedral solution of HWP$(48; 3, 4; 13, 10)$}

\noindent
Consider the nine cycles of $K_{48}$ defined as follows.
$$
\begin{array}{ll}
C_1=&\bigl{(}1, \ -\frac 1{\sqrt 2}(i+j), \ -\frac 1{\sqrt 2}(1+j)\bigl{)}\smallskip\\
C_2=&\bigl{(}1, \ \frac 1{2}(1-i+j-k), \ -\frac 1{\sqrt 2}(i+k)\bigl{)}\smallskip\\
C_3=&\bigl{(}1, \ \frac 1{\sqrt 2}(-i+j), \ \frac 1{2}(1-i-j-k)\bigl{)}\smallskip\\
C_4=&\bigl{(}1, \ \frac 1{2}(-1+i-j+k), \ \frac 1{\sqrt 2}(i-k)\bigl{)}\smallskip\\
C_5=&\bigl{(}1, \ \frac 12(-1-i+j+k), \ \frac 12(-1+i-j-k)\bigl{)}\smallskip\\
C_6=&\bigl{(}1, \ -\frac 1{2}(1+i+j+k), \ \frac 1{2}(-1+i-j+k), \ \frac 1{\sqrt 2}(1+j)\bigl{)}\smallskip\\
C_7=&\bigl{(}1, \ -\frac 1{\sqrt 2}(1+k), \ -k, \ \frac 1{2}(-1+i+j-k)\bigl{)}\smallskip\\
C_8=&\bigl{(}1, \ k, \ -1, \ -k\bigl{)}\smallskip\\
C_9=&\bigl{(}1, \ j, \ -1, \ -j\bigl{)}\\
\end{array}
$$

We note that for $5\leq i\leq 7$ the $G$-stabilizer of $C_i$ is $V(C_i)$ while all other $C_i$'s have trivial $G$-stabilizer.
Thus, using partial differences, one can check that $Orb_G(C_i)$ is a $\ell_i$-cycle decomposition of Cay$[G:\Omega_i]$ where $\ell_i$
is the length of $C_i$ and where the $\Omega_i$'s
are the symmetric subsets of $G$ listed below.

$$
\begin{array}{ll}
\Omega_1=&\{-\frac 1{\sqrt 2}(i+j), -\frac 1{\sqrt 2}(1+j), \frac 1{2}(1+i+j-k)\}^{\pm 1}\smallskip\\
\Omega_2=&\{\frac 1{2}(1-i+j-k), -\frac 1{\sqrt 2}(i+k), \frac 1{\sqrt 2}(1+i)\}^{\pm 1}\smallskip\\
\Omega_3=&\{\frac 1{\sqrt 2}(-i+j), \frac 1{2}(1-i-j-k), \frac 1{\sqrt 2}(j-k)\}^{\pm 1}\smallskip\\
\Omega_4=&\{\frac 1{2}(-1+i-j+k), \frac 1{\sqrt 2}(i-k), -\frac 1{\sqrt 2}(j+k)\}^{\pm 1}\smallskip\\
\Omega_5=&\{\frac 12(-1-i+j+k)\}^{\pm 1}\smallskip\\
\Omega_6=&\{k\}^{\pm 1}\smallskip\\
\Omega_7=&\{j\}^{\pm 1}\smallskip\\
\Omega_8=&\{-\frac 1{2}(1+i+j+k), i, \frac 1{\sqrt 2}(-1+i), \frac 1{\sqrt 2}(1+j)\}^{\pm 1}\smallskip\\
\Omega_9=&\{-\frac 1{\sqrt 2}(1+k), \frac 1{\sqrt 2}(1-k), \frac 1{2}(1-i+j+k), \frac 12(-1+i+j-k)\}^{\pm 1}\smallskip\\
\end{array}
$$
Now note that the $\Omega_i$'s partition $G\smallsetminus\{1, -1\}$.
This assures that \break $\displaystyle {\cal C}:= \bigcup_{i=1}^{9}Orb_{G}(C_i)$ is a $G$-regular
cycle-decomposition of $K_{48}-I$.

Now set $F_i=Orb_{S_i}(C_i)$ where
$S_i=\begin{cases}K & \mbox{for $1 \leq i \leq 4$;}\cr G & \mbox{for $5\leq i\leq 7$;}\cr L & \mbox{for $i=8, 9$}.\end{cases}$

Each $F_i$ is a 2-factor of $K_{48}$ with $Stab_{G}(F_i)=S_i$, hence $Orb_{G}(F_i)$ has length
3 or 1 or 4 according to whether $1\leq i\leq 4$ or $5 \leq i \leq 7$ or $i=8, 9$, respectively.

The cycles of $F_i$ are triangles or quadrangles according to whether or not $i\leq5$.
Thus, recalling that $\cal C$ is a cycle-decomposition of $K_{48}-I$, we conclude that
$\displaystyle {\cal F}:= \bigcup_{i=1}^{9}Orb_{G}(F_i)$ is a $G$-regular
2-factorization of $K_{48}-I$ with 13 triangle-factors and 10 quadrangle-factors, namely a
$G$-regular solution of HWP$(48; 3, 4; 13, 10)$.

\subsection{An octahedral solution of HWP$(48; 3, 4; 15, 8)$}

\noindent
Consider the seven cycles of $K_{48}$ defined as follows.
$$
\begin{array}{ll}
C_1=&\bigl{(}1, \ \frac 12(-1-i+j+k), \ \frac 1{\sqrt 2}(i+k)\bigl{)}\smallskip\\
C_2=&\bigl{(}1, \ -\frac 1{\sqrt 2}(i+j), \ -\frac 1{\sqrt 2}(1+j)\bigl{)}\smallskip\\
C_3=&\bigl{(}1, \ \frac 12(-1+i+j-k), \ \frac 12(1-i+j+k)\bigl{)}\smallskip\\
C_4=&\bigl{(}1, \ \frac 12(1+i+j+k), \ \frac 1{\sqrt 2}(1+j)\bigl{)}\smallskip\\
C_5=&\bigl{(}1, \ \frac 12(1-i+j-k), \ \frac 1{\sqrt 2}(i-k)\bigl{)}\smallskip\\
C_6=&\bigl{(}1, \ -j, \ k, \ -\frac 1{\sqrt 2}(1-k)\bigl{)}\smallskip\\
C_7=&\bigl{(}1, \ \frac 1{\sqrt 2}(i-j), \ \frac 12(-1-i+j-k), \ \frac 12(-1+i+j+k)\bigl{)}\\
\end{array}
$$

Here, every $C_i$ has trivial $G$-stabilizer. Thus, using partial differences,
one can check that $Orb_G(C_i)$ is a $\ell_i$-cycle decomposition of Cay$[G:\Omega_i]$ where $\ell_i$
is the length of $C_i$ and where the $\Omega_i$'s
are the symmetric subsets of $G$ listed below.

$$
\begin{array}{ll}
\Omega_1=&\{\frac 12(-1-i+j+k), \frac 1{\sqrt 2}(i+k), \frac 1{\sqrt 2}(-j+k)\}^{\pm 1}\smallskip\\
\Omega_2=&\{-\frac 1{\sqrt 2}(i+j), -\frac 1{\sqrt 2}(1+j), \frac 12(1+i+j-k)\}^{\pm 1}\smallskip\\
\Omega_3=&\{\frac 12(-1+i+j-k), \frac 12(1-i+j+k), \frac 12 (-1-i+j-k)\}^{\pm 1}\smallskip\\
\Omega_4=&\{\frac 12(1+i+j+k), \frac 1{\sqrt 2}(1+j), \frac 1{\sqrt 2}(1+i)\}^{\pm 1}\smallskip\\
\Omega_5=&\{\frac 12(1-i+j-k), \frac 1{\sqrt 2}(i-k), \frac 1{\sqrt 2}(j+k)\}^{\pm 1}\smallskip\\
\Omega_6=&\{-j, +i, \frac 1{\sqrt 2}(1-k), -\frac 1{\sqrt 2}(1-k)\}^{\pm 1}\smallskip\\
\Omega_7=&\{\frac 1{\sqrt 2}(i-j), -\frac 1{\sqrt 2}(1+i), +k, \frac 12(-1+i+j+k)\}^{\pm 1}\smallskip\\
\end{array}
$$

Now note that the $\Omega_i$'s partition $G\smallsetminus\{1, -1\}$.
This assures that \break $\displaystyle {\cal C}:= \bigcup_{i=1}^{8}Orb_{G}(C_i)$ is a $G$-regular
cycle-decomposition of $K_{48}-I$.

Set $F_i=Orb_{S_i}(C_i)$ where
$S_i=\begin{cases}K & \mbox{for $1 \leq i \leq 5$;}\cr L & \mbox{for $i=6, 7$}.\end{cases}$

Each $F_i$ is a 2-factor of $K_{48}$ with $Stab_{G}(F_i)=S_i$, hence $Orb_{G}(F_i)$ has length
3 or 4 according to whether $1\leq i\leq 5$ or $i=6,7$, respectively.
The cycles of $F_i$ are triangles or quadrangles according to whether or not $i\leq5$.
Thus, recalling that $\cal C$ is a cycle-decomposition of $K_{48}-I$, we conclude that
$\displaystyle {\cal F}:= \bigcup_{i=1}^{7}Orb_{G}(F_i)$ is a $G$-regular
2-factorization of $K_{48}-I$ with 15 triangle-factors and 8 quadrangle-factors, namely a
$G$-regular solution of HWP$(48; 3, 4; 15, 8)$.

\subsection{An octahedral solution of HWP$(48; 3, 4; 17, 6)$}

\noindent
Consider the ten cycles of $K_{48}$ defined as follows.
$$
\begin{array}{ll}
C_1=&\bigl{(}1,  \ -\frac 1{\sqrt 2}(1-k), \  -\frac 1{\sqrt 2}(i+k)\bigl{)}\smallskip\\
C_2=&\bigl{(}1, \ -\frac 1{\sqrt 2}(i+j), \ \frac 1{2}(-1+i+j+k)\bigl{)}\smallskip\\
C_3=&\bigl{(}1, \ \frac 12(1+i-j-k), \ -\frac 1{\sqrt 2}(1+j)\bigl{)}\smallskip\\
C_4=&\bigl{(}1, \ \frac 1{\sqrt 2}(-i+j), \ \frac 1{\sqrt 2}(-i+k)\bigl{)}\smallskip\\
C_5=&\bigl{(}1, \ \frac 12(1-i+j-k), \ \frac 1{\sqrt 2}(1-j)\bigl{)}\smallskip\\
C_6=&\bigl{(}1, \ \frac 12(-1-i+j+k), \ \frac 12(-1+i-j-k)\bigl{)}\smallskip\\
C_7=&\bigl{(}1, \ \frac 12(-1+i+j-k), \ \frac 12(-1-i-j+k)\bigl{)}\smallskip\\
C_8=&\bigl{(}1, \ k, \ -1, \ -k\bigl{)}\smallskip\\
C_9=&\bigl{(}1, \ j, \ -1, \ -j\bigl{)}\smallskip\\
C_{10}=&\bigl{(}1, \ \frac 1{\sqrt 2}(1+i), \ \frac 1{\sqrt 2}(1-i), \frac 1{2}(1-i-j+k)\bigl{)}\smallskip\\
\end{array}
$$

We note that for $6\leq i\leq 9$ the $G$-stabilizer of $C_i$ is $V(C_i)$ while all other $C_i$'s have trivial $G$-stabilizer.
Thus, using partial differences, one can check that $Orb_G(C_i)$ is a $\ell_i$-cycle decomposition of Cay$[G:\Omega_i]$ where $\ell_i$
is the length of $C_i$ and where the $\Omega_i$'s
are the symmetric subsets of $G$ listed below.

$$
\begin{array}{l}
\Omega_1=\{-\frac 1{\sqrt 2}(1-k), -\frac 1{\sqrt 2}(i+k), \frac 1{2}(-1-i+j-k)\}^{\pm 1}\smallskip\\
\Omega_2=\{-\frac 1{\sqrt 2}(i+j), \frac 1{2}(-1+i+j+k), \frac 1{\sqrt 2}(-1+i)\}^{\pm 1}\smallskip\\
\Omega_3=\{\frac 12(1+i-j-k), -\frac 1{\sqrt 2}(1+j), \frac 1{\sqrt 2}(j+k)\}^{\pm 1}\smallskip\\
\Omega_4=\{\frac 1{\sqrt 2}(-i+j), \frac 1{\sqrt 2}(-i+k), \frac 1{2}(1-i-j-k)\}^{\pm 1}\smallskip\\
\Omega_5=\{\frac 12(1-i+j-k), \frac 1{\sqrt 2}(1-j), \frac 1{\sqrt 2}(j-k)\}^{\pm 1}\smallskip\\
\Omega_6=\{\frac 12(-1-i+j+k)\}^{\pm 1}\smallskip\\
\Omega_7=\{\frac 12(-1+i+j-k)\}^{\pm 1}\smallskip\\
\Omega_8=\{k\}^{\pm 1}\smallskip\\
\Omega_9=\{j\}^{\pm 1}\smallskip\\
\Omega_{10}=\{\frac 1{\sqrt 2}(1+i), i, \frac 1{\sqrt 2}(1-k), \frac 12(1-i-j+k)\}^{\pm 1}\smallskip\\
\end{array}
$$

Now note that the $\Omega_i$'s partition $G\smallsetminus\{1, -1\}$.
This assures that \break $\displaystyle {\cal C}:= \bigcup_{i=1}^{8}Orb_{G}(C_i)$ is a $G$-regular
cycle-decomposition of $K_{48}-I$.

Set $F_i=Orb_{S_i}(C_i)$ where
$S_i=\begin{cases}K & \mbox{for $1 \leq i \leq 5$;}\cr G & \mbox{for $6 \leq i \leq 9$;}\cr L & \mbox{for $i=10$}.\end{cases}$

Each $F_i$ is a 2-factor of $K_{48}$ with $Stab_{G}(F_i)=S_i$, hence $Orb_{G}(F_i)$ has length
3 or 1 or 4 according to whether $1\leq i\leq 5$ or $6 \leq i \leq 9$ or $i=10$, respectively.
The cycles of $F_i$ are triangles or quadrangles according to whether or not $i\leq7$.
Thus, recalling that $\cal C$ is a cycle-decomposition of $K_{48}-I$, we conclude that
$\displaystyle {\cal F}:= \bigcup_{i=1}^{10}Orb_{G}(F_i)$ is a $G$-regular
2-factorization of $K_{48}-I$ with 17 triangle-factors and 6 quadrangle-factors, namely a
$G$-regular solution of HWP$(48; 3, 4; 17, 6)$.

\smallskip

\section {Dicyclic solutions of two Hamilton-Waterloo problems}

In this section $G$ will denote the dicyclic group of order 24 which is usually denoted by $Q_{24}$.
Thus $G$ has the following presentation:

$$
G=\langle a, b\, |\, a^{12}=1,\, b^2=a^6,\, b^{-1}a
b=a^{-1}\rangle
$$

\noindent Note that the elements of $G$ can be written in the form $a^ib^j$ with $0\leq i\leq 11$ and $j=0, 1$. The group $G$
has a unique involution which is $a^6$ and we will need to consider the following subgroups of $G$:
\begin{itemize}
\item $H=\langle b\rangle=\{1, b, a^6, a^6b\}$;
\item $K=\langle a^2\rangle=\{1, a^2, a^4, a^6, a^8, a^{10}\}$;
\item $L=\langle a^2b, a^3\rangle=\{1, a^3, a^6, a^9, a^2b, a^8b, a^5b,
a^{11}b\}$.
\end{itemize}

In the following the complete graph $K_{24}$ and the cocktail party graph $K_{24}-I$
will be seen as the Cayley graphs Cay$[G:G\setminus\{1\}]$ and \break Cay$[G:G\setminus\{1,a^6\}]$, respectively.

\subsection{A dicyclic solution of HWP$(24; 3, 4; 7, 4)$}

\noindent
Consider the four cycles of $K_{24}$ defined as follows.
$$
\begin{array}{ll}
C_1=&\bigl{(}1, \ a^3b, \ a^5\bigl{)}\smallskip\\
C_2=&\bigl{(}1, \ a^{10}, \ a^7b\bigl{)}\smallskip\\
C_3=&\bigl{(}1, \ a^4, \ a^8\bigl{)}\smallskip\\
C_4=&\bigl{(}1, \ b, \ a^3b, \ a\bigl{)}\smallskip\\
\end{array}
$$

We note that the $G$-stabilizer of $C_3$ is $V(C_3)$ while all other $C_i$'s have trivial $G$-stabilizer.
Thus, using partial differences, one can check that $Orb_G(C_i)$ is a $\ell_i$-cycle decomposition of Cay$[G:\Omega_i]$ where $\ell_i$
is the length of $C_i$ and where the $\Omega_i$'s
are the symmetric subsets of $G$ listed below.

$$
\begin{array}{ll}
\Omega_1=&\{a^3b, a^5, a^2b\}^{\pm 1}\smallskip\\
\Omega_2=&\{a^2, ab, a^5b\}^{\pm 1}\smallskip\\
\Omega_3=&\{a^4\}^{\pm 1}\smallskip\\
\Omega_4=&\{b, a^3, a^4b, a\}^{\pm 1}\smallskip\\
\end{array}
$$

Now note that the $\Omega_i$'s partition $G\smallsetminus\{1, a^6\}$.
This assures that \break $\displaystyle {\cal C}:= \bigcup_{i=1}^{4}Orb_G(C_i)$ is a $G$-regular
cycle-decomposition of $K_{24}-I$.

Now set $F_i=Orb_{S_i}(C_i)$ where
$S_i=\begin{cases}L & \mbox{for $i=1, 2$;}\cr G & \mbox{for $i=3$;}\cr K & \mbox{for $i=4$}.\end{cases}$

Each $F_i$ is a 2-factor of $K_{24}$ with $Stab_G(F_i)=S_i$, hence $Orb_G(F_i)$ has length
3 or 1 or 4 according to whether $i=1, 2$ or $i=3$ or $i=4$, respectively.

The cycles of $F_i$ are triangles or quadrangles according to whether or not $i\leq3$.
Thus, recalling that $\cal C$ is a cycle-decomposition of $K_{48}-I$, we conclude that
$\displaystyle {\cal F}:= \bigcup_{i=1}^{4}Orb_{G}(F_i)$ is a $G$-regular
2-factorization of $K_{24}-I$ with 7 triangle-factors and 4 quadrangle-factors, namely a
$G$-regular solution of HWP$(24; 3, 4; 7, 4)$.

\subsection{A dicyclic solution of HWP$(24; 3, 4; 9, 2)$}

\noindent
Consider the four cycles of $K_{24}$ defined as follows.

$$
\begin{array}{ll}
C_1=&\bigl{(}1, b, a^6, a^6b\bigl{)}\smallskip\\
C_2=&\bigl{(}1, a^4b, a^6, a^{10}b\bigl{)}\smallskip\\
C_3=&\bigl{(}1, a^4, a^7b\bigl{)}\smallskip\\
C_4=&\bigl{(}1, a^3b, a^8b\bigl{)}\smallskip\\
C_5=&\bigl{(}a^4, a^7, a^5\bigl{)}\smallskip\\
\end{array}
$$

We note that for $i=1, 2$ the $G$-stabilizer of $C_i$ is $V(C_i)$ while all other $C_i$'s have trivial $G$-stabilizer.
Thus, using partial differences, one can check that $Orb_G(C_i)$ is a $\ell_i$-cycle decomposition of Cay$[G:\Omega_i]$
where $\ell_i$ is the length of $C_i$ and where the $\Omega_i$'s are the symmetric subsets of $G$ listed below.

$$
\begin{array}{ll}
\Omega_1=&\{b\}^{\pm 1}\smallskip\\
\Omega_2=&\{a^4b\}^{\pm 1}\smallskip\\
\Omega_3=&\{a^4, ab, a^5b\}^{\pm 1}\smallskip\\
\Omega_4=&\{a^3b, a^2b, a^5\}^{\pm 1}\smallskip\\
\Omega_5=&\{a^1, a^2, a^3\}^{\pm 1}\smallskip\\
\end{array}
$$

Also here the $\Omega_i$'s partition $G\smallsetminus\{1, a^6\}$, hence
$\displaystyle {\cal C}:= \bigcup_{i=1}^{5}Orb_{G}(C_i)$ is a $G$-regular
cycle-decomposition of $K_{24}-I$. Now set:

$$F_1=Orb_{G}(C_1),\quad\quad\quad\quad\quad\quad\quad F_2=Orb_{G}(C_2),$$
$$F_3=Orb_L(C_3),\quad\quad F_4=Orb_H(C_4) \ \cup \ Orb_H(C_5).$$

Each $F_i$ is a 2-factor of $K_{24}$ and we have
$$Stab_G(F_1)=Stab_G(F_2)=G; \quad Stab_G(F_3)=L; \quad Stab_G(F_4)=H$$
so that the lengths of the $G$-orbits of $F_1$, \dots, $F_4$ are 1, 1, 3 and 6, respectively.
The cycles of $F_i$ are triangles or quadrangles according to whether or not $i\geq3$.
Thus, recalling that $\cal C$ is a cycle-decomposition of $K_{48}-I$, we conclude that
$\displaystyle {\cal F}:= \bigcup_{i=1}^{5}Orb_{G}(F_i)$ is a $G$-regular
2-factorization of $K_{24}-I$ with 9 triangle-factors and 2 quadrangle-factors, namely a
$G$-regular solution of HWP$(24; 3, 4; 9, 2)$.

\smallskip

\section {A special linear solution of HWP$(24; 3, 4; 5, 6)$}

In this section $G$ will denote the 2-dimensional special linear group over $\mathbb{Z}_3$, usually denoted by $SL_2(3)$,
namely the group of $2\times2$ matrices with elements in $\mathbb{Z}_3$ and determinant one.
The only involution of $G$ is $2E$ where $E$ is the identity matrix of $G$.
The 2-Sylow subgroup $Q$ of $G$,
isomorphic to the group of quaternions, is the following:
$$Q=\biggl{\{}  \begin{bmatrix}1&0\\0&1\\ \end{bmatrix}, \begin{bmatrix}2&0\\0&2\\ \end{bmatrix},
\begin{bmatrix}1&1\\1&2\\ \end{bmatrix}, \begin{bmatrix}2&2\\2&1\\ \end{bmatrix}, \begin{bmatrix}0&2\\1&0\\ \end{bmatrix},
\begin{bmatrix}0&1\\2&0\\ \end{bmatrix}, \begin{bmatrix}1&2\\2&2\\ \end{bmatrix}, \begin{bmatrix}2&1\\1&1\\ \end{bmatrix}\biggl{\}}.$$

\noindent We will also need to consider the subgroup $H$ of $G$ of order 6 generated by the matrix $\begin{bmatrix}0&1\\2&1\\ \end{bmatrix}$.
Hence we have:
$$H=\biggl{\{}  \begin{bmatrix}1&0\\0&1\\ \end{bmatrix}, \begin{bmatrix}0&1\\2&1\\ \end{bmatrix},
\begin{bmatrix}2&1\\2&0\\ \end{bmatrix}, \begin{bmatrix}2&0\\0&2\\ \end{bmatrix}, \begin{bmatrix}0&2\\1&2\\ \end{bmatrix},
\begin{bmatrix}1&2\\1&0\\ \end{bmatrix}\biggl{\}}.$$

The use of the special linear group $G$ was crucial in \cite{B2004} to get a Steiner triple system of any order $v=144n+25$
with an automorphism group acting sharply transitively an all but one point.
Here $G$ will be used to get a $G$-regular solution of the last Hamilton-Waterloo problem left open in \cite{DQS}.
In the following the complete graph $K_{24}$ and the cocktail party graph $K_{24}-I$
will be seen as the Cayley graphs Cay$[G:G\setminus\{E\}]$ and \break Cay$[G:G\setminus\{E,-E\}]$, respectively.

Consider the six cycles of $K_{24}$ defined as follows.

$$C_1=\biggl{(}\begin{bmatrix}1&0\\0&1\\
\end{bmatrix}, \
\begin{bmatrix}2&0\\
2&2\\
\end{bmatrix}, \
\begin{bmatrix}1&2\\
1&0\
\end{bmatrix}\biggl{)}$$

$$C_2=\biggl{(}\begin{bmatrix}1&0\\
0&1\\
\end{bmatrix}, \
\begin{bmatrix}0&2\\
1&2\\
\end{bmatrix}, \
\begin{bmatrix}2&1\\
2&0\\
\end{bmatrix}\biggl{)}$$

$$C_3=\biggl{(}\begin{bmatrix}
1&0\\
0&1\\
\end{bmatrix}, \
\begin{bmatrix}
0&1\\
2&2\\
\end{bmatrix}, \
\begin{bmatrix}
2&2\\
1&0\\
\end{bmatrix}\biggl{)}$$

$$C_4=\biggl{(}\begin{bmatrix}
1&0\\
0&1\\
\end{bmatrix}, \
\begin{bmatrix}
1&0\\
2&1\\
\end{bmatrix}, \
\begin{bmatrix}
2&2\\
0&2\\
\end{bmatrix}, \
\begin{bmatrix}
1&1\\
0&1\\
\end{bmatrix}\biggl{)}$$

$$C_5=\biggl{(}\begin{bmatrix}
1&0\\
0&1\\
\end{bmatrix}, \
\begin{bmatrix}
1&1\\
1&2\\
\end{bmatrix}, \
\begin{bmatrix}
2&0\\
0&2\\
\end{bmatrix}, \
\begin{bmatrix}
2&2\\
2&1\\
\end{bmatrix}\biggl{)}$$

$$C_6=\biggl{(}\begin{bmatrix}
1&0\\
0&1\\
\end{bmatrix}, \
\begin{bmatrix}
2&1\\
1&1\\
\end{bmatrix}, \
\begin{bmatrix}
2&1\\
0&2\\
\end{bmatrix}, \
\begin{bmatrix}
1&1\\
0&1\\
\end{bmatrix}\biggl{)}$$

\medskip\noindent
Here the $G$-stabilizer of $C_i$ is trivial for $i=1, 6$ while it coincides with $V(C_i)$ for $2\leq i\leq 5$.
Thus, using partial differences, one can check that $Orb_G(C_i)$ is a $\ell_i$-cycle decomposition of Cay$[G:\Omega_i]$ where $\ell_i$
is the length of $C_i$ and where the $\Omega_i$'s
are the symmetric subsets of $G$ listed below.

$$\Omega_1=\biggl{\{}\begin{bmatrix}
2&0\\
2&2\\
\end{bmatrix}, \
\begin{bmatrix}
1&2\\
1&0\\
\end{bmatrix}, \
\begin{bmatrix}
0&2\\
1&1\\
\end{bmatrix}\biggl{\}}^{\pm1}$$

$$\Omega_2=\biggl{\{}\begin{bmatrix}
0&2\\
1&2\\
\end{bmatrix}\biggl{\}}^{\pm1}\quad\quad
\Omega_3=\biggl{\{}\begin{bmatrix}
0&1\\
2&2\\
\end{bmatrix}\biggl{\}}^{\pm1}$$
$$\Omega_4=\biggl{\{}\begin{bmatrix}
0&1\\
2&0\\
\end{bmatrix}\biggl{\}}^{\pm1}\quad\quad
\Omega_5=\biggl{\{}\begin{bmatrix}
1&1\\
1&2\\
\end{bmatrix}\biggl{\}}^{\pm1}$$

$$\Omega_6=\biggl{\{}\begin{bmatrix}
2&1\\
1&1\\
\end{bmatrix}, \
\begin{bmatrix}
1&0\\
2&1\\
\end{bmatrix}, \
\begin{bmatrix}
2&2\\
0&2\\
\end{bmatrix},
\begin{bmatrix}
1&1\\
0&1\\
\end{bmatrix}\biggl{\}}^{\pm1}$$

Once again we see that the $\Omega_i$'s partition $G\smallsetminus\{E,2E\}$, therefore\break
$\displaystyle {\cal C}:= \bigcup_{i=1}^{4}Orb_G(C_i)$ is a $G$-regular
cycle-decomposition of $K_{24}-I$.

\noindent\medskip
Now set $F_i=Orb_{S_i}(C_i)$ with
$S_i=\begin{cases}Q & \mbox{for $i=1$;}\cr G & \mbox{for $2\leq i\leq 5$;}\cr H & \mbox{for $i=6$}.\end{cases}$

Each $F_i$ is a 2-factor of $K_{24}$ and we have
$Stab_{G}(F_i)=S_i$ so that the lengths of the $G$-orbits of $F_1$, \dots, $F_6$ are 3, 1, 1, 1, 1 and 4, respectively.

The cycles of $F_i$ have length 3 or 4 according to whether or not $i\leq 3$ or not.
Thus, recalling that $\cal C$ is a cycle-decomposition of $K_{24}-I$, we conclude that
$\displaystyle {\cal F}:= \bigcup_{i=1}^{6}Orb_{G}(F_i)$ is a $G$-regular
2-factorization of $K_{24}-I$ with 5 triangle-factors and 6 quadrangle-factors, namely a
$G$-regular solution of HWP$(24; 3, 4; 5, 5)$.


\begin{thebibliography}{99}

\bibitem{BurBon} S. Bonvicini and M. Buratti, {\it Sharply vertex transitive factorizations of Cayley graphs}, in preparation.

\bibitem{BBRT} S. Bonvicini, M. Buratti, G. Rinaldi and T. Traetta,
{\it Some progress on $1$-rotational Steiner triple systems},
Des. Codes Cryptogr. {\bf62} (2012), 63-78.

\bibitem{B2001} M. Buratti,
{\it $1$-rotational Steiner triple systems
over arbitrary groups}, J. Combin. Des. {\bf 9} (2001), 215-226.

\bibitem{B} M. Buratti,
{\it Rotational $k$-cycle systems of order $v<3k$; another proof of the existence of odd cycle systems},
J. Combin. Des. {\bf11} (2003), 433--441.

\bibitem{B2004} M. Buratti,
{\it Cycle decompositions with a sharply vertex transitive automorphism group},
Matematiche (Catania) {\bf59} (2004), 91--105.

\bibitem{BR} M. Buratti and G. Rinaldi,
{\it On sharply vertex transitive $2$-factorizations of the complete graph},
J. Combin. Theory Ser. A {\bf111} (2005), 245--256.

\bibitem{DQS} P. Danziger, G. Quattrocchi, B. Stevens, {\it The Hamilton-Waterloo problem for cycle Sizes $3$ and $4$},
J. Combin. Des. {\bf17} (2009), 342--352.


\end{thebibliography}
\end{document}